\documentclass[12pt,leqno]{amsart}

\usepackage[latin2]{inputenc}
\usepackage{amsmath}
\usepackage{amsthm}
\usepackage{amssymb}
\usepackage[mathscr]{eucal}
\usepackage{enumerate}
\usepackage{amsfonts}
\usepackage{geometry,amssymb}
 \textwidth=16cm

\newtheorem{theorem}{Theorem}
\newtheorem{lemma}{Lemma}

\begin{document}
\baselineskip=15pt
\title[Large Oscillation of Fourier Transforms]
{Oscillation of Fourier Transforms and Markov-Bernstein Inequalities}

\author[Sz. Gy. R\'ev\'esz, N. N. Reyes, G. A. M. Velasco]
{Szil\'ard Gy. R\'ev\'esz, Noli N. Reyes and  Gino Angelo M. Velasco}

\address{Szil\'ard Gy. R\'ev\'esz \newline\indent A. R\'enyi Institute of Mathematics \newline \indent Hungarian
Academy of Sciences, \newline \indent Budapest, P.O.B. 127, 1364
\newline \indent Hungary} \email{revesz@renyi.hu}

\address{Noli N. Reyes
\newline\indent Department of Mathematics
\newline\indent College of Science
\newline\indent University of the Philippines
\newline\indent Quezon City, 1101, Philippines}
\email{noli@math.upd.edu.ph}


\begin{abstract}
\baselineskip=15pt Under certain conditions on an integrable function $ P
$ having a real-valued Fourier transform $ \hat{ P }$ and such that $ P(0)
= 0 $, we obtain an estimate which describes the oscillation of $ \hat{P}
$ in the interval $ [ -C \| P' \|_{ \infty} /  \| P \|_{\infty} , C \| P'
\|_{ \infty} /  \| P \|_{\infty}] $, where $C$ is an absolute constant,
independent of $P$. Given $ \lambda > 0 $ and an integrable function $
\phi $ with a non-negative Fourier transform, this estimate allows us to
construct a finite linear combination $ P_{ \lambda} $ of the translates $
\phi ( \cdot + k \lambda ) , \ k \in {\bf Z}$ such that $  \| P_{ \lambda}
'\|_{ \infty } >  c \| P_{ \lambda} \|_{ \infty } / \lambda  $ with
another absolute constant $c>0$. In particular, our construction proves
sharpness of an inequality of Mhaskar for Gaussian networks.
\end{abstract}



\thanks{The second author acknowledges the support of the Natural Sciences
Research Institute of the University of the Philippines.}

\thanks{The third author was supported in part through the Hungarian-French
Scientific and Technological Governmental Cooperation, Project \# F-10/04,
the Hungarian-Spanish Scientific and Technological Governmental
Cooperation, Project \# E-38/04 and by the Hungarian National Foundation
for Scientific Research, Project \#s T-049301, T-049693 and K-61908.}

\maketitle


\section{Introduction}
\setlength\baselineskip{15pt}

The original A.A. Markov inequality  states that
$   || P'||_{ L^{ \infty }(I) }  \leq n^{2} || P ||_{ L^{ \infty }(I) } $
for any algebraic polynomial $P$ of degree $n$. Here, $ I= [-1,1]  $.
This inequality becomes an equality if $P$ is the Chebyshev polynomial
$  P(x) = \cos nt $ where $ x = \cos t  $. The reader may find the details of this
in  page 40 of \cite{lorentz}.

Upper estimates of the derivative norm by that of the function itself are usually termed
Markov-Bernstein inequalities.
There is an extensive literature on such inequalities, which
play an important role in inverse theorems, where
smoothness of a function is deduced from rates of convergence of polynomial approximations.
For an excellent survey on Markov-Bernstein and related inequalities, the reader may consult
the book \cite{borwein} of P. Borwein and T. Erd\'elyi.

By imposing additional assumptions on the zeros of the polynomials, one can obtain
estimates which give lower estimates for  the norm of a
derivative in terms of the norm of the function. These results are usually
termed inverse Markov-Bernstein inequalities or Tur\'an type inequalities.
 For instance, Tur\'an \cite{turan} proved that
\[ || P'||_{ L^{ \infty }(I) }   \geq  \frac{ \sqrt{n}}{ 6 } || P ||_{ L^{ \infty }(I) }  \]
 for any polynomial $P$ of degree $n$, provided that all of its zeros lie in the interval
$I = [-1,1]$. We also refer the reader to a valuable paper of Er\"od \cite{erod}.

There is an upsurge of interest in such estimates, with a number of
recent results dealing with the topic (\cite{erdelyi2}, \cite{levenberg},
\cite{szilard}, \cite{zhou}).
For instance, in  \cite{zhou}, Zhou showed that if $ 0 <r \leq q \leq \infty $ and
$ 1 \geq 1/r - 1/q$, then
\[ || P'||_{ L^{ r }(I) }   \geq C  n^{\alpha} || P ||_{ L^{ q }(I) }  \]
for every polynomial $P$ whose zeros lie in the interval $I$. Here,
$ \alpha = \frac{1 }{2 } - \frac{ 1}{2r } + \frac{1 }{ 2q} $.

More related to our work are results of  Erd\'elyi and Nevai \cite{erdelyi}
where they obtained
\[  \lim_{ n \rightarrow \infty} \frac{ ||p_{n}'||_{X} }{ ||p_{n} ||_{Y} } =  \infty \]
for sequences of polynomials $ p_{n}$ whose zeros satisfy cetain conditions .

Markov-Bernstein inequalities have also been obtained for other classes of functions such as
Gaussian networks. For instance, in \cite{mhaskar}, Mhaskar showed that
for some constant $c$, $ ||g' ||_{p} \leq c m || g ||_{p} $ for any function
$g$ defined on the real line of the form
\[ g(x)  \ = \ \sum_{k=1}^{N} a_{k} \exp ( - ( x-x_{k} )^{2}  ),  \]
where $ | x_{j} - x_{k}  | \geq 1/m $ for $  j \neq k $, and $ \log N =
{\mathcal O} (m^{2}) $.

One of our  goals in this note is to show that under certain conditions on
an integrable function $ P :{\bf R} \longrightarrow {\bf R} $ having a
real-valued Fourier transform $ \hat{ P }$ with $ P(0) = 0 $,
\begin{equation}
r \geq  C  \frac{\| P'\|_{ \infty } }{\| P\|_{ \infty }} \  \
\Longrightarrow  \ \  \int_{ -r}^{r } ( \hat{P })_{ \pm } \ \geq
\frac{\sqrt{2\pi}}{4}\ ||  P ||_{ \infty }  .
\label{eqn:oscillation}
\end{equation}
Here, we can take $C= 8^{3} / \pi $.  This estimate not only tells us that
$ \hat{P} $ will have a zero in the interval $ [-r,r] $, but also provides
an effective estimate on how it oscillates in the interval.

For a fixed function $ \phi $, let
\begin{equation}\label{eqn:Endef}
E_n(\lambda):=\left\{ \sum_{k= -n}^{n } b_{k} \phi (x + \lambda k)
~:~ b_k\in {\bf R},~k=-n,\dots,-1,0,1,\dots,n \right\}.
\end{equation}
The estimate in (\ref{eqn:oscillation}) allows us to construct
  $ P_{ \lambda } \in E_n( \lambda ) $ for each $ \lambda > 0$
  and for sufficiently large positive integers $n$ (depending on $ \lambda$)
such that $  \| P_{ \lambda} '\|_{ \infty }  > c
  \| P_{ \lambda} \|_{ \infty }  / \lambda $ with some absolute constant $c>0$.
In particular, our construction proves sharpness of the above-mentioned
inequality of Mhaskar \cite{mhaskar} for Gaussian networks.

\section{Notations and preliminaries}

For any integrable function $f$ on the real line, we write for its
Fourier transform
\[  \hat{ f }( \omega ) \ = \ \frac{1}{ \sqrt{ 2 \pi}}
  \int_{ {\bf R}} f(x) e^{ -i \omega x} \ dx \ . \]
Given a real number $x$, its positive and negative parts are
$x_{+} = \max \{ x,0 \} $ and $  x_{-} = \max \{-x,0 \} $
respectively.

We will write $h$ for the Fej\'er kernel, that is
\[ h( x ) \ := \ \frac{1}{\sqrt{ 2 \pi } }
 \left( \frac{ \sin x/2 }{  x/2 }   \right)^{2} \ . \]
Its Fourier transform is given by
\[  \hat{ h } ( \omega )  \ = \
  \max \{ 1 - | \omega |, \ 0    \} \ . \]

For the rest of the paper we fix an auxiliary function $H$. We could use
any even, $ 2 \pi $-periodic, and e.g. twice continuously differentiable
function \\ $ H: {\bf R} \longrightarrow {\bf R} $, not identically one, such
that $ H(x) = 1  $ if $ |x| \leq \pi/2$. The special constants and values
in the following choice are not relevant, only some order is essential.
Nevertheless, for definiteness and more explicit calculation we take e.g.
\begin{equation}
H(x)= \left\{
  \begin{array}{ll}
    1, \qquad\qquad & \mbox{if} \ |x| \leq \pi /2;  \\
    \sin^2x, & \mbox{if } \pi/2 < |x| \leq \pi .
\end{array}
\right.  \label{eqn:Hdef}
\end{equation}
Then $H$ has the Fourier cosine series development
\[ H(x) \ = \ \sum_{k=0}^{\infty } a_{k} \cos kx  \]
where $a_k$ are the Fourier cosine coefficients  of $H$. Although
precise values are not needed here, a calculation leads to
$a_0=3/4$, $a_1=4/(3\pi)$, $a_2=\frac{-1}{4}$ and
\begin{equation}
a_k= \frac{-4\sin\frac{k\pi}{2}}{\pi k(k^2-4)} =
\left\{
\begin{array}{ll}
    0, & k \mbox{ even;} \\
    \frac{-4}{\pi k(k^2-4)}, &  k \equiv 1 \mbox{ mod } 4; \\
    \frac{4}{\pi k(k^2-4)}, & k\equiv 3 \mbox{ mod } 4.
\end{array}
\right. ~~~~ \hbox{for}~k \geq 3, ~ k\in {\bf N}.
\label{eqn:akdef}
\end{equation}
It is immediate that $|a_k| \leq k^{-2} $ for all $k \in {\bf N}$;
moreover, a direct calculation yields
\begin{equation}
\sum_{k=1}^{\infty}  |a_k| = 1 +\frac{5}{3\pi} = 1.530516...<1.6
\qquad \textrm{and}\qquad \sum_{k=1}^{\infty} a_k^2 = \frac{9}{8}.
\label{eqn:akabssum} \end{equation}

\section{Oscillation of Fourier transforms}

\begin{lemma}\label{l:oscillation}
Let $ P: {\bf R} \longrightarrow {\bf R} $ be bounded,
differentiable, and integrable such that  $ \hat{ P } $ is
real-valued. Suppose $ P(0) = 0$ and let
\begin{equation}
r\ > \ \frac{8^{3} \ ||P'||_{\infty }}{\pi\ || P||_{\infty}} \ ,
\label{eqn:lambda}
\end{equation}
then
\begin{equation}
\frac{4}{\sqrt{ 2 \pi}}  \int_{ -r}^{r } ( \hat{P })_{\pm} \ \geq
\ || P ||_{ \infty }
\end{equation}
\end{lemma}

{\bf Proof of Lemma \ref{l:oscillation}}: There is nothing to
prove if $ || P' ||_{ \infty} = \infty $. Hence, we assume $ || P'
||_{ \infty} < \infty $. Fix $ r $ satisfying (\ref{eqn:lambda})
and define
\[ f(x) \ = \ P \star h_{ r } (x) \ = \  \frac{1}{ \sqrt{ 2 \pi }}
 \int_{ {\bf R}}    P( x-t ) h_{ r } (t) dt     \]
where $ h_{ r } (t) = r h ( r t) $.
Since $  ( 2 \pi )^{ -1/2} \int_{ {\bf R}} h_{ r } =1 $,
for any real number $x$,
\begin{equation}
 f(x) - P(x) = S(x) + L(x)
\label{eqn:split}
\end{equation}
where
\[ S(x) =  \frac{1}{ \sqrt{ 2 \pi }} \int_{  |t| < \delta } ( P(x-t ) - P(x)  )
h_{ r } (t) \ dt \ , \]
\[  L(x) =  \frac{1}{ \sqrt{ 2 \pi }} \int_{  |t| \geq \delta } ( P(x-t ) - P(x)  )
h_{ r } (t) \ dt \  \]
and $ \delta > 0 $ is chosen  such that $ 8 \delta q = 1  $ with
$ q = || P'||_{ \infty } / || P ||_{ \infty } $.
Combining the inequalities
\[ | S(x)| \ \leq \ \delta \ || P' ||_{ \infty } \ = \
 \frac{ || P ||_{ \infty }}{ 8}
\ \ \mbox{ and } \ \ | L(x) | \ \leq \
\frac{ 8 \ || P||_{ \infty } }{ \pi r \delta } \ < \
\frac{|| P||_{ \infty }}{ 8 }  \]
with (\ref{eqn:split}), we obtain for any real number $x$,
\begin{equation}
 |f(x) - P(x)  | < || P||_{ \infty } / 4 \ .
\label{eqn:main}
\end{equation}

Since $f$ and $ \hat{f} $ are both integrable,  the inversion
formula for the Fourier transform shows that
\[  \sqrt{ 2 \pi } || f ||_{ \infty } \ \leq \  \int_{ {\bf R}} | \hat{ f }  |
   \ = \ \int_{ {\bf R}} \left( \hat{ f }  + 2  ( \hat{ f } )_{-} \right)
 \  =  \ \sqrt{ 2 \pi }  f(0) + 2 \int_{ {\bf R}}( \hat{ f } )_{-} .    \]
Applying (\ref{eqn:main}) with $ x= 0$ and noting that $ P(0) = 0
$, we conclude that
\[
|| f ||_{ \infty } \ \leq \ \frac{1}{4}|| P ||_{ \infty } +
  \frac{2}{ \sqrt{2 \pi } } \int_{ {\bf R}}( \hat{ f } )_{-} \ .
\]
Making use once more of (\ref{eqn:main}) and the last inequality gives
\[
|| P ||_{ \infty } \ \leq \ || f ||_{ \infty } + \frac{1}{4}|| P
||_{ \infty } \ \leq \ \frac{1}{2}|| P ||_{ \infty } + \frac{2}{
\sqrt{ 2 \pi } } \int_{ {\bf R}}( \hat{ f } )_{-}
\]
and therefore
 \[  || P ||_{ \infty } \  \leq \
  \frac{4}{\sqrt{ 2 \pi}} \int_{ {\bf R}}( \hat{ f } )_{-} \ . \]

 Finally, we observe that $ \hat{f} ( \omega ) = \hat{P} ( \omega ) \hat{h} ( r^{-1} \omega )$,
  $ 0 \leq  \hat{h} \leq 1 $ and  $ \hat{h} = 0 $ outside $  [-1,1]$. These
  imply that   $   ( \hat{ f } )_{-} = 0   $  outside $ [-r,r] $
  and $ ( \hat{ f } )_{-} \leq ( \hat{ P})_{- } \  $.
  Therefore
  \[   || P ||_{ \infty } \ \leq \
   \frac{4}{\sqrt{ 2 \pi}} \int_{ -r}^{r } ( \hat{ P } )_{-} \ .  \]
  A similar argument leads to the same inequality for $ ( \hat{ P } )_{+} \ $.
 $ \Box $


\section{Construction of sums of translates with \\ large oscillation}


\begin{theorem}
Let $ \phi : {\bf R} \longrightarrow {\bf R}$ be an even, continuous,
integrable function such that $ \phi (0) =1 $. In addition, suppose that
its Fourier transform $ \hat{ \phi } $ is nonnegative, integrable and
analytic on $ {\bf R} $. Given $ \lambda > 0 $, then there exist a
positive integer $n$ and $P\in E_n(\lambda)$, with $E_n(\lambda)$ defined
in (\ref{eqn:Endef}), such that
\[ \frac{|| P '||_{\infty }}{|| P ||_{\infty }} \ \geq
 \frac{ C}{ \lambda} . \] Here, we could take $ C = \pi^{2}/2^{10} $.
\end{theorem}

{\bf Proof:} For each positive integer $n$ and for each real number $x$,
we define
\begin{equation}\label{eqn:pndef}
P_{n}(x) = 2 A_{n} \phi (x) \ + \ \sum_{ k=1}^{n} a_{k} ( \phi ( x
+ \lambda k)  + \phi ( x - \lambda k) )
\end{equation}
and also
\begin{equation}
P_{ \infty } (x) \ := \ \lim_{ n \rightarrow \infty} P_{n} (x) =2
A_{\infty}(\lambda)\phi (x) \ + \ \sum_{ k=1}^{\infty} a_{k} (
\phi ( x + \lambda k)  + \phi ( x - \lambda k) ),
\label{eqn:pinfty}
\end{equation}
where the coefficients $a_k$ are the Fourier cosine coefficients
of $H$ in (\ref{eqn:akdef}), and
\begin{equation}\label{eqn:Andef}
A_{n} := A_n(\lambda):=- \sum_{k=1}^{n} a_{k} \phi ( \lambda k),
\qquad A_{\infty} := A_\infty(\lambda):=- \sum_{k=1}^{\infty}
a_{k} \phi (\lambda k).
\end{equation}

We start with showing that $ P_{ \infty } $ is not identically zero.
\begin{lemma}\label{l:Pnotvanish}
Under the assumptions of Theorem 1, we have $ || P_{ \infty } ||_{
\infty } > 0 $.
\end{lemma}

{\bf Proof of lemma \ref{l:Pnotvanish}:} For each  $ \omega \in
{\bf R} $ and $n\in{\bf N}$ we define
\begin{equation}
T_{n} ( \omega ) \ := \ \sum_{k=1}^{n} a_{k} ( \cos ( k \lambda
\omega ) - \phi ( \lambda k ) ) \label{eqn:trig}
\end{equation}
and also
\begin{equation}
T_{ \infty }( \omega ) \ := \ \lim_{n \rightarrow \infty} T_{n} (
\omega ) = \sum_{k=1}^{\infty} a_{k} ( \cos ( k \lambda \omega ) -
\phi ( \lambda k ) ). \label{eqn:triginfty}
\end{equation}
Thus, $  \hat{ P}_{ \infty } ( \omega ) \ = \hat{\phi}(\omega) 2 T_{
\infty }( \omega ) = 2 \hat{ \phi } ( \omega ) ( H( \lambda \omega ) - F(
\lambda ) ) $, where
\begin{equation}\label{eqn:Fdef}
F( \lambda) := \sum_{ k=0}^{\infty}  a_{k} \phi(\lambda k ) =
a_0-A_{\infty}(\lambda)~
\end{equation}
is a uniformly convergent sum of bounded functions of $\lambda$.
By the Fourier inversion formula, $ P_{ \infty } \equiv 0 $ if and
only if $ \hat{P}_{ \infty} \equiv 0 $. Thus, it suffices to show
that for any given $ \lambda
> 0$, $ \hat{ \phi } ( \omega ) ( H( \lambda \omega ) - F( \lambda
) ) $ does not vanish identically.

Note that for any $ \lambda > 0$, $ H ( \lambda \omega) \neq 1 $
for  $ \ \omega \in {\mathcal I} = (  \frac{ \pi }{ 2 \lambda  } ,
\frac{ 3\pi }{ 2 \lambda  } )
     +  (2 \pi /\lambda)  {\bf Z} $ while
    $ H ( \lambda \omega) = 1 $ for
    $ \ \omega \in {\mathcal J} = [ - \frac{ \pi }{ 2 \lambda  } ,
    \frac{ \pi }{ 2 \lambda  } ]
     +  (2 \pi /\lambda) {\bf Z} $.
Therefore, if  $ F( \lambda) =1 $, then $ F( \lambda ) \neq H(
\lambda \omega) $ for $ \omega  \in {\mathcal I} $, while if $ F(
\lambda)  \neq 1 $, then $ F( \lambda ) \neq H( \lambda \omega) $
for $ \omega  \in {\mathcal J} $. In any case, $ H( \lambda \omega
) - F( \lambda )  \neq 0 $ for $ \omega $ in a union of non-empty
open intervals. If $ \hat{P}_{ \infty} \equiv 0$, then $ \hat{
\phi } $ would have to be zero on these intervals, which is
impossible since $ \hat{ \phi }  $ is assumed to be  analytic on $
{\bf R}$. This completes the proof of lemma 2. $ \Box $

To finish the proof of the theorem it suffices to show the next
assertion.

\begin{lemma}\label{claim:punchline} If a positive integer  $n$ is chosen
such that $ 20 \sum_{ k > n } | a_{k}| < || P_{ \infty}  ||_{ \infty} $,
then
\[ \frac{ || P_{n} '||_{ \infty }}{ || P_{n} ||_{ \infty }} \ \geq
\ \frac{ \pi^{2} }{  2^{10} \lambda  }. \]
\end{lemma}
{\bf Proof of Lemma \ref{claim:punchline}:} Recall $\hat{P}_{n} (
\omega ) = 2 \hat{ \phi} ( \omega )  T_{n} ( \omega )  $ with
$T_{n} $ in (\ref{eqn:trig}). We also define $ \Delta_{n} ( \omega
) = T_{n} ( \omega ) - H ( \lambda \omega ) + F( \lambda) $ for $
\omega\in {\bf R}$, with $F(\lambda)$ in (\ref{eqn:Fdef}).

Meanwhile, in view of the assumptions that $ \hat{ \phi } \geq 0 $
and $ \hat{ \phi } \in L^{1} $, the inversion formula for the
Fourier transform shows that $ ||\phi||_{\infty} = \phi(0) = 1$.
With this in mind, we obtain for every positive integer $n$
\begin{equation}
|| \Delta_{n}  ||_{ \infty } \leq \ 2 \sum_{ k > n } | a_{k}|~,
     \qquad \mbox{and} \qquad
||P_{\infty} - P_{n} ||_{ \infty } \leq \ 4 \sum_{k > n } |a_{k}|.
\label{eqn:error}
\end{equation}

Suppose $  0 < r  \leq \pi / ( 2 \lambda ) $. Then $ H( \lambda \omega  ) =1  $
for $  | \omega | \leq r$.
Therefore, if $ F( \lambda )  \geq 1 $, then
\[ \int_{ -r}^{r } ( \hat{ P}_{n} )_{ +} \ = \ 2  \int_{ -r}^{r } \hat{ \phi
}~
  ( 1 - F( \lambda )+ \Delta_n )_{ +}
  \ \leq \ 4  \sqrt{ 2 \pi }  \sum_{ k > n } | a_{k}| . \]
   Here, we've again made use of the conditions $  \hat{ \phi } \geq 0$ and
   $ \phi (0 ) =1 $.
   Similarly, if $ F( \lambda )  < 1 $, we also obtain
 \[ \int_{ -r}^{r } ( \hat{ P}_{n} )_{- }  \ \leq \ 4
   \sqrt{ 2 \pi }  \sum_{ k > n } | a_{k}| . \]
Thus, we've shown that  if $  0 < r  \leq \pi / ( 2 \lambda ) $, then
 for each positive integer $n$,
 \begin{equation}
 \min \left( \int_{ -r}^{r } ( \hat{ P}_{n} )_{- }, \int_{ -r}^{r }
 ( \hat{ P}_{n} )_{+ }  \right)
   \ \leq \ 4 \sqrt{ 2 \pi }   \sum_{ k > n } | a_{k}|.
   \label{eqn:summary} \end{equation}
On the other hand, lemma \ref{l:oscillation} together with the
second inequality in (\ref{eqn:error}) asserts that if $r > (8^{3}
/ \pi) || P_{n}' ||_{ \infty } /  || P_{n}||_{ \infty }$, then
\begin{equation}\label{eqn:lower-bound}
\frac{ 4}{ \sqrt{ 2 \pi} } \min \left( \int_{ -r}^{r } ( \hat{
P}_{n} )_{- }~, \int_{ -r}^{r } ( \hat{ P}_{n} )_{+ }  \right)
\geq || P_{n}||_{ \infty}  \geq  || P_{\infty}||_{ \infty} - 4
\sum_{ k
> n } | a_{k}|.
\end{equation}
Combining (\ref{eqn:summary}) and (\ref{eqn:lower-bound}) we
conclude that if $ (8^{3} / \pi) || P_{n}'||_{ \infty} / || P_{n}
||_{ \infty} < \pi /(2 \lambda)$, then $ || P_{ \infty } ||_{
\infty }- 4\sum\limits_{ k > n } | a_{k}| \leq 16 \sum\limits_{ k
> n } | a_{k}|$ and therefore
  $ || P_{ \infty } ||_{ \infty } \leq 20 \sum\limits_{ k > n } | a_{k}| $.
This proves the lemma which gives the conclusion of the theorem.
$\Box$

 \section{Application to Gaussian networks}

Our goal in this section is to prove sharpness of an inequality of
Mhaskar (mentioned in the introduction of this paper) for Gaussian
networks. We shall apply Theorem 1 (in particular, lemma 1 in the
proof) with $ \phi (x ) = \exp ( -x^{2} ) $. In this section
$E_n(\lambda)$ is defined according to (\ref{eqn:Endef}) with our
above given Gaussian $\phi$.

The following theorem is the main result of this section.

\begin{theorem}
Let $ n \in {\bf N} $ and $ \lambda \in (0,1) $ satisfy
\begin{equation}\label{eqn:Nzerocond}
n > N_0:= C_0  \lambda \exp\left( \frac{\pi^2}{2\lambda^2}
\right)\qquad \qquad \left( C_0:=\frac{1280}{3\pi}\right).
\end{equation}
Then there exists $P\in E_n(\lambda)$, such that
\[
\frac{|| P '||_{ \infty }}{|| P ||_{ \infty }}
 \geq    \frac{ \pi^{2} }{ 2^{10} \lambda  }.
 \]
\end{theorem}

{\bf Remark.} Note $\log N_0 = O (1/\lambda^2)$, in complete
agreement with the above mentioned result of H. N. Mhaskar. Thus
the result proves sharpness of the result in \cite{mhaskar} for an
arithmetic progression of shifts $x_k:=\lambda k$ with separation
$1/m=\lambda$.
\[
\]
We retain the function $H$ from (\ref{eqn:Hdef}) and its Fourier coefficients
$a_k$ in (\ref{eqn:akdef}) also in this section. With these Fourier
coefficients $a_k$ and for each $ \lambda > 0 $ and $x \in {\bf R}$,
$P_{\infty}(\lambda,x)$ will again be as in (\ref{eqn:pinfty}) with $
A_{\infty} ( \lambda ) $ defined in (\ref{eqn:Andef}). However, in
contrast to the proof of Theorem 1, $ \lambda  $ is no longer fixed.

As we are dealing with the Gaussian function
$\phi(x):=\exp(-x^2)$, a number of properties are immediate.

First of all, the fact that $ \phi $ is even and decreasing on
$[0, \infty )$ implies that for each $ \lambda > 0 $ and for any
real number $x$,
\begin{equation} \ \sum_{ k \in {\bf Z}}  \phi ( k \lambda - x  )
\ \leq \ 1 + \frac{1}{\lambda } \int_{ {\bf R}} \phi = 1 +
\frac{\sqrt{\pi}}{\lambda}. \label{eqn:riemann}
\end{equation}
Indeed, all values of $\phi ( k \lambda - x  )$ can be replaced by
the $\int$ over the interval of length $\lambda$ from $k \lambda - x$
towards $0$, except perhaps the function value at the (single, if $x\ne
\pm\lambda/2$) point which is closest to $0$ (and thus is estimated by 1).

Also, the Fourier transform of $\phi$ is given by $ \hat{ \phi} ( \omega )
= (1/\sqrt{ 2 }) \exp ( - \omega ^{2} / 4) $. Keeping only the term with
maximal absolute value, we easily obtain
\begin{equation}\label{eqn:phihatest}
\sum_{l=-\infty}^{\infty} \left|\hat{\phi}  \left( \frac{\omega +
2\pi l} {\lambda} \right)   \right|^{2} \geq  \hat{\phi}^{2}
\left(\frac{\pi}{\lambda}\right) \qquad\qquad \left( \forall
~\omega \in {\bf R}\right).
\end{equation}

\begin{lemma} For the function (\ref{eqn:pinfty}) we have
\begin{equation}  |P_{\infty}(
\lambda, x)| \leq \frac{24}{1+x^2} \qquad \qquad \left( x\in {\bf
R}\right), \label{eqn:decay-of-p}
\end{equation} uniformly for all $\lambda \in (0,1)$.
\label{lemma:unidecay}
\end{lemma}
{\bf Proof of Lemma \ref{lemma:unidecay}:} Using $\phi(\lambda k)\leq 1$
and (\ref{eqn:akabssum}) we obtain
$$
|A_\infty(\lambda)| \leq \sum_{k=1}^\infty |a_k|  \leq 1.6~.
$$
As $\max\limits_{\bf R}(1+x^2)\phi(x)= \max\limits_{[0,\infty)}
(1+t)e^{-t}=1$, we get
\begin{equation}\label{Aestim}
|2A_\infty(\lambda)\phi(x)| \leq \frac{3.2}{1+x^2}.
\end{equation}
It follows that we indeed have
\begin{equation} |P_{\infty}( \lambda, x)| \ \leq \frac{3.2}{1+x^2} +
\  \sum_{k \in {\bf Z}\setminus 0} |a_k| \phi ( x - \lambda k),
\label{eqn:decroiss}
\end{equation}
where $ a_{k} = a_{-k} $ if $ k <0$. As $\|\phi\|_{\infty} = \phi(0)= 1$,
in case $|x|\le 2$ this immediately leads to $|P_{\infty}( \lambda, x)|
\leq {3.2}/(1+x^2) + 3.2 < 20/(1+x^2)$, hence (\ref{eqn:decay-of-p}).

Because the right hand side of (\ref{eqn:decroiss}) is even, it
remains to take $x>2$.

Now let $\mathcal A$ be the set of all nonzero integers $k$ such that
$|x-\lambda k | < x /2 $. Observe that for $k\in\mathcal A~$, $\lambda |k|
\geq x/2 $ and thus $|k| \geq x/(2\lambda)$, which gives by $|a_k|\leq
1/k^2$, also $|a_k| \leq 4\lambda^2/x^2 \leq 5 \lambda^{2}/(1+x^2)$ for
$x>2$. Therefore, taking into account also (\ref{eqn:riemann}) and $x>2$,
we are led to
\begin{equation}
\sum_{ k \in {\mathcal A}} |a_k| \phi ( x - \lambda k) \leq \frac{5
\lambda^2}{1+x^2} \left( 1 + \frac{\sqrt{\pi}}{\lambda} \right) \leq
\frac{5\lambda^2+5 {\sqrt{\pi}}{\lambda}}{1+x^2} . \label{eqn:nearx}
\end{equation} On the other hand, in view of (\ref{eqn:akabssum}) and
\[ \max\limits_{[2,\infty)} ( 1+x^{2}  ) \phi(x/2)=
\max\limits_{[4,\infty)}(1+t)e^{-t/4}=5/e  , \]
we have
\begin{equation}
\sum_{ k \not\in {\mathcal A} } |a_k|\phi ( x - \lambda k) \leq \
\phi\left(\frac{x}{2}\right) 2 \sum_{k=1}^{\infty} |a_k| \leq
\frac{10}{e(1+x^{2})}  \sum_{k=1}^{\infty} |a_k| <
\frac{6}{1+x^2}. \label{eqn:farx}
\end{equation}
Recalling $0<\lambda<1$ a combination of (\ref{eqn:decroiss}),
(\ref{eqn:nearx}) and (\ref{eqn:farx}) gives the result of the
lemma. $ \Box$
\[  \]

We  shall also make use of an explicit lower bound for the
$L^{2}$-norm of $  P_{ \infty } ( \lambda , \cdot ) $ in terms of
the $ l^{2}$-norm of its coefficients. Actually, a more general
phenomenon can be observed here.
\begin{lemma}  Let $ \lambda > 0 $ be fixed and $c_k\in{\bf C}$
$(k\in{\bf Z})$ be arbitrary coefficients satisfying $\sum_{k\in{\bf Z}}
|c_k|^{2} < \infty$, i.e., $(c_k)\in\ell_2({\bf Z})$. Consider the
function $f(\lambda,x):=\sum_{k=-\infty}^{\infty} c_k \phi(x-\lambda k)$.
We then have
\begin{equation}\label{eqn:independence}
  || f (\lambda, \cdot )  ||^{2}_{2}
    \ \geq \  \mu ( \lambda ) \sum_{ k =-\infty}^{ \infty} |c_{k}|^{2}
\end{equation}
where
\begin{equation}
\mu ( \lambda ) : = \ \frac{ 2 \pi }{ \lambda } \inf_{ \omega \in {\bf R}}
\sum_{ l \in {\bf Z}} \left| \hat{\phi} \left( \frac{\omega + 2 \pi
l}{\lambda } \right) \right|^{2}. \label{eqn:riesz-bound}
\end{equation}
\label{lemma:independence}
\end{lemma}
{\bf Proof of Lemma
\ref{lemma:independence}:} First of all, for a fixed $ \lambda > 0$, the
series defining $ f:=f (\lambda , \cdot ) $ converges in $ L^{2} ( {\bf
R}) $. To see this, we consider its sequence
 $ f_{n} ( \lambda , x ) = \sum_{|k| \leq n} c_k \phi(x-\lambda k)$
 of partial sums. The Fourier transform of $ f_{n}:=f_n (\lambda ,
\cdot ) $ is given by  $   \hat{f}_{n} ( \lambda , t) = \hat{ \phi
}(t) \sum_{|k| \leq n } c_{k} e^{ - ik \lambda t}$. Applying
Plancherel's theorem to
$ || f_{n}(\lambda , \cdot ) - f_{m}(\lambda , \cdot )||^{2}_{2} $ and
writing the resulting integral as a sum of integrals over the
intervals $  [ 2 \pi \lambda^{-1}l, 2 \pi \lambda^{-1}(l+1)  ] , \
l \in {\bf Z} $, we obtain
\[ || f_{n}(\lambda , \cdot ) - f_{m}(\lambda , \cdot )  ||^{2}_{2} \ = \
 \frac{1}{ \lambda } \sum_{ l = - \infty}^{ \infty }
  \int_{ 0}^{ 2 \pi }  \left|  \hat{ \phi }
   \left( \frac{ \omega + 2 \pi l}{ \lambda } \right)
  \sum_{ m < |k| \leq n }   c_k e^{ ik \omega} \right|^{2}
    d \omega       \]
    for $m < n $.
    The rapid decay of $ \hat{ \phi }$  assures the
    finiteness of
    \[ M(\lambda) \ := \ \frac{2 \pi}{ \lambda } \
    \sup\limits_{ \omega \in {\bf R}} \
     \sum_{ l = - \infty }^{  \infty}      \left|  \hat{ \phi }
   \left( \frac{ \omega + 2 \pi l}{ \lambda } \right)
  \right|^{2}   \]
  and therefore by Parseval's theorem,
  \[ || f_{n}(\lambda , \cdot ) - f_{m}(\lambda , \cdot )  ||^{2}_{2}
   \ \leq \  M_{ \lambda } \sum\limits_{m< |k| \leq n }
    | c_{k}|^{2}  \ \longrightarrow 0 \]
    as $ n > m \longrightarrow  \infty $.
   This proves convergence in $ L^{2} $ of the series defining
   $ f (\lambda , \cdot ) $.

A similar argument furnishes the  conclusion of the lemma except that we
take the infimum $\mu(\lambda)$ (as defined in (\ref{eqn:riesz-bound})),
 instead of the supremum $M(\lambda)$ above. $  \Box$


{\bf Proof of Theorem 2:} First of all, we estimate $ ||
P_{\infty}( \lambda,  \cdot ) ||_{ \infty } $ from below by $ ||
P_{\infty}( \lambda,  \cdot ) ||_{ 2 } $. In view of lemma
\ref{lemma:unidecay} we have
\begin{equation}
|P_{\infty}( \lambda, x )| \leq \frac{C}{|x|} \qquad
 (\textrm{with}~~C=12)\label{eqn:decay}
\end{equation}
for all real numbers $x\ne 0$ and for each $\lambda >0$.

Now let the parameter $\sigma$ be chosen so that
$$
\sigma := \frac{4C^{2}}{|| P_{\infty}(\lambda, \cdot ) ||_{ 2 }^2} .
$$
Note that $P_{\infty}$ does not vanish identically, hence
$ \sigma >0$.
We write $ || P_{\infty}( \lambda, \cdot ) ||_{ 2 }^{2} $ as a sum
of integrals over $ [- \sigma , \sigma ] $ and over $ {\bf R} \setminus
 [- \sigma , \sigma ]$. Estimating trivially in
 $[- \sigma, \sigma ]$ and applying (\ref{eqn:decay})
to the second integral yields
 \[  || P_{\infty}( \lambda,  \cdot ) ||_{ 2 }^{2} \leq
 2 \sigma || P_{\infty}( \lambda,  \cdot ) ||_{ \infty }^{2} +
  2 C^{2} \sigma^{-1} . \]
Thus, a short calculation with the chosen value of $ \sigma$ leads for
each $ \lambda> 0$
\begin{equation}
   || P_{\infty}( \lambda,  \cdot ) ||_{ 2 }^{2} \ \leq \ 4C
   || P_{\infty}( \lambda,  \cdot ) ||_{ \infty }.
\label{eqn:infty-two}
\end{equation}

To evaluate $\|P_{\infty}(\lambda,\cdot)\|_2$ we note
$P_{\infty}(\lambda,x)= \sum_{ k \in {\bf Z} } \alpha_{k} \phi(x-
k \lambda)$, where $\alpha_{k} = a_{|k|}$ if $ k \neq 0$, and $
\alpha_{0} = 2 A_{\infty}(\lambda)$, with $A_{\infty}(\lambda)$
defined in (\ref{eqn:Andef}). For this function we clearly have
$\sum_{k\in{\bf Z}} |\alpha_k|^2 \geq 2 \sum_{k=1}^{\infty}
|a_k|^2=9/4$ in view of (\ref{eqn:akabssum}).

Meanwhile, we consider the function $ \mu(\lambda) $ defined in
(\ref{eqn:riesz-bound}). Recalling (\ref{eqn:phihatest}) and the explicit
form of $\widehat{\phi}$ provides for each $ \lambda > 0$ the estimate
\[
\mu( \lambda )  \geq  \frac{\pi}{ \lambda } \exp \left( -
\frac{\pi^2}{2\lambda^2} \right).
\]
Combining this with lemma \ref{lemma:independence} and
(\ref{eqn:infty-two}) we obtain
\begin{equation}
|| P_{\infty}( \lambda,  \cdot ) ||_{ \infty } \geq \frac{\pi}{4C\lambda}
\exp \left( - \frac{\pi^2}{2\lambda^2} \right) \sum_{k=-\infty}^{\infty}
|\alpha_k|^2 = \frac{9\pi}{16C\lambda} \exp \left( -
\frac{\pi^2}{2\lambda^2} \right)~. \label{eqn:l2}
\end{equation}

Now recalling $|a_k|\leq 1/k^2$ we obtain $ \sum_{ k>n} |a_{k}| < 1/n$ for
each positive integer $n$. Recalling also $C=12$, this and (\ref{eqn:l2})
yields that whenever (\ref{eqn:Nzerocond}) holds, then
\[
20 \sum_{k>n} |a_{k}| < 20/n < 20/N_0 =  \frac{3\pi }{64 \lambda}
\exp\left(-\frac{\pi^2}{2\lambda^2}\right) < || P_{\infty}(
\lambda,  \cdot ) ||_{ \infty }.
\]
Therefore, an application of lemma \ref{claim:punchline} concludes
the proof of Theorem 2. $ \Box $

\end{document}